\documentclass[11pt]{article}
\catcode`\@=11 \catcode`\@=12

\usepackage{amsmath,amscd}
\usepackage{amssymb}
\usepackage{amsthm}
\usepackage{oldgerm}

\usepackage[pdftex,colorlinks,urlcolor=blue,pdfstartview=FitH]{hyperref}

\newcommand{\Rm}{\mathbb{R}}

\newcommand{\mL}{\mathcal{L}}

\newcommand{\mC}{\ensuremath{\mathcal{C}}}
\newcommand{\mR}{\ensuremath{\mathcal{R}}}

\newcommand{\mX}{\ensuremath{\mathcal{X}}}

\newcommand{\Nm}{\ensuremath{\mathbb{N}}}
\newcommand{\Zm}{\ensuremath{\mathbb{Z}}}

\newcommand{\mM}{\ensuremath{\mathcal{M}}}

\newcommand{\mI}{\ensuremath{\mathcal{I}}}

\newcommand{\mG}{\ensuremath{\mathcal{G}}}

\newtheorem{lem}{Lemma}
\newtheorem{thm}{Theorem}

\newtheorem{cor}[lem]{Corollary}
\newtheorem{prop}[lem]{Proposition}
\newtheorem{defn}[lem]{Definition}

\def\proof {\noindent{\sc{Proof. }}}
\def\qed {\mbox{}\hfill {\small \fbox{}} \\}
\def\lto{\longrightarrow}
\def\lmto{\longmapsto}

\def\leq{\leqslant}
\def\geq{\geqslant}

\title{On the number of Mather measures of Lagrangian systems
}
\author{Patrick  Bernard}
\date{ November 2009 \footnote{This is version 2, the first version was submitted 
in  July 2008}}
\begin{document}

\maketitle
\begin{small}

 Patrick Bernard
 \footnote{membre de l'IUF}\\
CEREMADE, CNRS\\
Pl. du Mar\'echal de Lattre de Tassigny\\
75775 Paris Cedex 16\\
France\\
\texttt{patrick.bernard@ceremade.dauphine.fr}\\

\begin{center}
-----
\end{center}

Abstract:
In 1996, Ricardo 
Ricardo Ma\~n\'e
discovered that Mather measures are in fact the 
minimizers  of a "universal" infinite dimensional
linear programming problem.
This fundamental result has many applications, one of the most 
important is to the  estimates of the generic number of 
Mather measures. 
Ma\~n\'e obtained the first estimation of that sort by using
finite dimensional approximations.
Recently, we were able with Gonzalo Contreras to use
this method of finite dimensional approximation 
in order to solve a conjecture of John Mather concerning
the generic number of Mather measures for families
of Lagrangian systems.
In the present paper we obtain finer results 
in that direction by
applying directly some  classical 
tools of  convex analysis to the infinite dimensional problem.
We  use   a notion of countably rectifiable sets of finite
codimension in  Banach (and Frechet)   spaces  which may deserve
 independent interest.

\begin{center}
-----
\end{center}

R\'esum\'e:
En 1996, Ricardo Ma\~n\'e a d\'ecouvert que les mesures de Mather peuvent
 \^etre obtenues comme
solutions d'un probl\`eme variationnel convexe "universel" de dimenion
infinie.
Ce r\'esultat fondamental a de nombreuses applications,
il permet par exemple d'estimer le nombre de mesures
de Mather des syst\`emes g\'en\'eriques.
Ma\~n\'e a obtenu la premi\`ere estimation de ce type 
en utilisant une approximation par des probl\`emes  variationnels
de dimension finie.
Nous avons r\'ecemment utilis\'e cette m\'ethode avec Gonzalo Contreras 
pour  r\'esoudre une conjecture de John
Mather sur le nombre g\'en\'erique de mesures minimisantes
dans les familles de syst\`emes Lagrangians.
Dans le pr\'esent article, on obtient des r\'esultat plus fins dans cette direction
en  appliquant  directement
au probl\`eme  de dimension infinie des m\'ethodes classiques de 
l'analyse convexe.  On \'etudie  pour ceci une nouvelle  notion 
d'ensembles rectifiables de codimension finie  dans les espaces de Banach
 (ou de Frechet)  qui est peut-\^etre int\'eressante en elle m\^eme.
\end{small}

\newpage

\section{Introduction}
Let $M$ be a compact connected manifold without boundary.
We want to study the dynamical system on $TM$
generated by a Tonelli Lagrangian 
$$
L:TM\lto \Rm.
$$
By \textit{Tonelli Lagrangian} we  mean a
  $C^2$ function $L:TM\lto \Rm$
such that, for each $x\in M$, the function
$v\lmto L(x,v)$ is superlinear and convex with positive definite Hessian.
Note then that the superlinearity is uniform with respect to $x$, see \cite{Fbook}, 
section 3.2.
To each Tonelli Lagrangian is associated
a complete $C^1$ flow $\psi^t$ on $TM$, with the property
that a curve $(x(t),v(t))$ is a trajectory of 
$\psi^t$ if and only if 
(i) $v(t)=\dot x(t)$  and (ii) the curve $x(t)$ solves the 
Euler-Lagrange equation
$$
\frac{d}{dt}(\partial_vL(x(t),\dot x(t)))=\partial_xL(x(t),\dot x(t)).
$$
We call this flow the Euler-Lagrange flow associated to $L$.
A standard example is the mechanical case where 
$$
L(x,v)=\frac{1}{2}\|v\|^2_x -V(x),
$$
the associated Euler Lagrange equation is just the Newton equation
$$
\ddot x(t)=-\nabla V (x(t)).
$$
Variational methods offer interesting tools to investigate the orbits
of the Euler-Lagrange flow.
We recall the following well-known fact: 
A $C^1$ curve $x(t)$
(or more generally an absolutely continuous curve $x(t)$)
 satisfies the Euler-Lagrange equation
on an open interval $I$ if and only if, for each $t_0\in I$,
there exists $\epsilon>0$ such that 
$[t_0-\epsilon,t_0+\epsilon]\subset I$
and such that 
$$\int_{t_0-\epsilon}^{t_0+\epsilon} L(x(t),\dot x(t))dt <
\int _{t_0-\epsilon}^{t_0+\epsilon} L(\gamma(t),\dot \gamma(t)) dt
$$
for all $C^1$ curves $\gamma:[t_0-\epsilon,t_0+\epsilon]\lto M$ 
different from $x$ and satisfying the boundary conditions
 $\gamma(t_0-\epsilon)=x(t_0-\epsilon)$
and $\gamma(t_0+\epsilon)=x(t_0+\epsilon)$.

One of the standard applications of variational
methods is to the existence of periodic orbits.
This can be done as follows. Fix a positive real number $T$
and a homology class $w\in H_1(M,\Zm)$.
Let $W^{1,1}(T,w)$ be the set of 
absolutely continuous curves $x:\Rm\lto M$
which are $T$-periodic and   have homology $w$ (when seen as closed 
loops on $M$).
It is a classical result that the action functional
$$
W^{1,1}(T,w)\ni \gamma\lmto \int_0^T L(\gamma(t),\dot \gamma(t))dt
$$
has a minimum, and that the minimizing curves are $C^2$
and solve the Euler-Lagrange equation.
This is a way to prove the existence of many periodic orbits
of the Euler-Lagrange flow.

John Mather had the idea to apply variational methods to
measures instead of curves.
Let $\mI(L)$ be the set of compactly supported Borel probability
measures on $TM$ which are invariant under the Euler-Lagrange flow.
Note that, if $x(t)$ is a $T$-periodic solution of the Euler-Lagrange
equation, then we associate to it an invariant measure $\mu$
characterized by the property that
$$
\int_{TM} f(x,v) d\mu(x,v)= \frac{1}{T}\int_0^T f(x(t),\dot x(t))dt
$$
for each continuous and bounded   function $f$ on $TM$.
In this case, we see that the action 
of the curve $x$ is just $T\int_{TM} Ld\mu$.
This suggests to take  $\int_{TM} Ld\mu$
 as the definition of the action of 
a probability measure.
A Mather measure is then  defined as a minimizer of the action  on
$\mI(L)$.
John Mather proved  in \cite{Mather:91}
that Mather measures exist, and moreover that they are
supported on a Lipschitz graph.
More precisely, there exists a Lipschitz vectorfield
$Y(x)$ on $M$
(which depends on the Lagrangian $L$, but not on the Mather measure)
 such that all the Mather measures of $L$
are supported on the graph of $Y$.

In the mechanical case where $L=\|v\|^2/2-V(x)$, 
 the Mather measures are just the invariant measures associated 
 to the rest points maximizing $V$
 (and therefore one can take $Y\equiv 0$ in this case).
 So we have not gained much insight in the dynamics of the 
 Euler-Lagrange flow of these Lagrangians at that point.
 A trick due to John Mather yet allows to obtain further information 
 from his theory. Recall first that, if $\omega$
 is a closed form on $M$, that we see as a function on $TM$
 linear in each fiber, then the Tonelli Lagrangian 
 $$\tilde L(x,v)=L(x,v)+\omega_x\cdot v
 $$
 generates the same Euler-Lagrange flow as $L$.
 This can be seen easily by considering the variational
 characterization of the Euler-Lagrange equation.
The remark of Mather is that,
 although $L$ and $\tilde L$
 generate the same flow, they 
 do not have the same Mather measures.
 Actually, the Mather measures of $\tilde L=L+\omega$
 depend 
  only on the cohomology of $\omega$
 in $H^1(M,\Rm)$. By definition, they are invariant measures
 of the flow of $L$. If the cohomology group $H^1(M,\Rm)$
 is not trivial, this construction allows to find
 non-trivial measures supported on Lipschitz graphs for
 mechanical Lagrangians.
 In order to simplify the notations for the sequel,
 we associate once and for all to each cohomology
 class $c\in H^1(M,\Rm)$ a smooth closed one-form
 $c_x\cdot v$ which has cohomology $c$. We assume that 
 the form depends linearly on the class.
  We denote by $L+c$
 the Lagrangian $L(x,v)+c_x\cdot v$.
 Finally, we denote by $\mM(L)$ the set of Mather measures of $L$.
 It is a convex subset of $\mI(L)$.

The Mather measures are very important objects in themselves,
 but they
also appear as building blocks for more elaborate constructions
of orbits of Lagrangian systems as initiated in \cite{Mather:93}
(see also \cite{fourier,JAMS,CY,CY2}).
It is  useful for these constructions to be able to 
understand the dimension of $\mM(L)$.
An important  result  was obtained by Ricardo Ma\~n\'e
in 1996, see \cite{Mane:96}:

\begin{thm}\label{Mane}
Let $L$ be a Tonelli Lagrangian, and let $\sigma^{\infty}(L)$
be the set of those potentials $V\in C^{\infty}(M)$
such that the Tonelli Lagrangian $L-V$ has more than
 one Mather measure.
The set $\sigma^{\infty}(L)$ is a meager set in the sense of Baire
category. It means that it is contained in the   union 
of countably many nowhere dense closed sets.
\end{thm}

When applied in the mechanical case, $L=\|v\|^2/2$
 this theorem states that 
generic smooth functions on $M$  have only one maximum,
which is of course not a new result.
More interesting situations appear by considering 
modified kinetic energies $L=\|v\|^2/2+\omega_x\cdot v$.
For  applications, however, it is necessary to treat simultaneously
all the sets $\mM(L-V+c)$, $c\in H^1(M,\Rm)$.
We were recently able with Gonzalo  Contreras to extend the result 
of Ma\~n\'e  in that direction, see  \cite{BC}.
These results imply Ma\~n\'e Theorem as well as
the following:

\begin{thm}\label{BC}
Let $L$ be a Tonelli Lagrangian.
Let $\Sigma_k^{\infty}(L)$ be the set of potentials 
$V\in C^{\infty}(M)$
such that, for some $c\in H^1(M,\Rm)$,
$dim(\mM(L-V+c))\geq k$.
Then, if $k>b_1$ (where $b_1=\dim H^1(M,\Rm)$), the set 
$\Sigma_k^{\infty}(L)$ is
Baire meager in $C^{\infty}(M)$.
\end{thm}

In the sequel, we shall fix $p\in \{2,3,\ldots ,\infty\}$
and consider potentials in $C^p(M)$
instead of $C^{\infty}(M)$ 
(note that the case $p=\infty$ is still included).
 We call 
$\Sigma_k^{p}(L)$ the set of potentials 
$V\in C^{p}(M)$
such that
$$
\max_{c\in H^1(M,\Rm)}\dim(\mM(L-V+c))\geq k.
$$
The results we have recalled above 
remain valid, with the same proof, in this context.
They state that it is exceptional in
the sense of Baire category to have too many Mather measures.
Now in the  separable Frechet space $C^p(M,\Rm)$,
 there are plenty of other
notions of small sets, which are at least as relevant as the
Baire category, and are more in the spirit of 
having  measure zero
(although there is not a single   way to define sets of measure
zero on infinite dimensional spaces).
Good introductions to these notions and to the literature concerning them
are \cite{BL}, \cite{Ph:78} and \cite{Cs:99}.
In dynamical systems, the most popular notion is prevalence,
that we now define:

A subset $A$ of a  Frechet space $B$ is said Haar-null
if there exists a compactly supported Borel probability measure $m$ on $B$
such that $m(A+x)=0$ for each $x\in B$.
This concept was first introduced by Christensen in the separable case,
 see \cite{Cr:72} or \cite{BL}.
It was  used as a description of the smallness
of the sets of non G\^ateau differentiability of Lipschitz functions
on separable Frechet spaces.
It was then rediscovered in the context of dynamical systems,
where the name \textit{prevalence} appeared,
see \cite{prev}. A prevalent set is the complement
of a Haar-null set.
It is proved in  \cite{prev}  that 
some versions of the Thom Transversality
Theorem still hold in the sense of prevalence.

Another notion is that of  Aronzsajn-null
sets, or equivalently of Gaussian-null sets, see \cite{Cs:99} or \cite{BL}.
They can be defined as those sets which have zero
measure for all Gaussian measure, see \cite{Ph:78}, \cite{Cs:99}
and Section \ref{small} 
for more details.
 The complement of an Aronzsajn-null set is prevalent.

One can  wonder whether the smallness results discussed
above concerning the dimension of $\mM(L)$
still hold for these notions of small sets.
Since these notions have first been introduced
to deal with non-differentiability points of Lipschitz or 
convex functions, and since 
the proof of the genericity results recalled above 
boils down to abstract convex analysis, it is not very
surprizing that the answer is positive.
It is  implied by the following stronger statement
expressed in terms of a notion of countably rectifiable sets 
that will be defined in Section \ref{sectrect}:

\begin{thm}\label{main}
Let $L$ be a Tonelli Lagrangian and let $p\in \{2,3,\ldots,\infty\}$ and $ k>b_1$
be given, where $b_1$ is the first Betti number of M (the dimension 
of   $H^1(M,\Rm)$).
The set  $\Sigma^p_k(L)$ is  countably rectifiable  of codimension
$k-b_1$ in $C^p(M)$. As a consequence, for each $k>b_1$, the set
$\Sigma_k^p(L)$ is meager in the sense of Baire and Aronszajn-null.
 Its complement is prevalent.
\end{thm}

We will give  in Section \ref{proof} a more general result which implies Theorem
\ref{main} and many similar statements.
 For example,
the set $\sigma^p_k(L)$ of those potentials $V\in C^p(M)$
such that $\dim(\mM(L-V))\geq k$ is countably rectifiable 
of codimension $k$. This is a refined version of Ma\~n\'e's result.

Here is an example of a new application:
In perturbation theory, one often fixes a Lagrangian $L$
and a potential $V$ and studies the dynamics generated
by $L-\epsilon V$, for $\epsilon$ small enough.
The following result is then useful:

\begin{cor}\label{epsilon}
Let $L$ be a Tonelli Lagrangian and let  $p\in \{2,3,\ldots,\infty\}$
be given.
Let $A_k^p(L)$ be the set of potentials $V\in C^p(M)$
such that
$$\sup_{\epsilon >0, c\in H^1(M,\Rm)}\dim(\mM(L-\epsilon V+c))\geq k.
$$
If $k>1+b_1$ ($b_1$ is the dimension of $H^1(M,\Rm)$), then $A_k^p(L)$ is 
countably rectifiable of codimension $k-(1+b_1)$ in $C^p(M)$.
As a consequence, the set $A^p_k(L)$ is Baire-meager and Aronszajn-null.
\end{cor}

This Corollary is proved in Section \ref{proof}.
Rectifiable sets of finite codimension in Banach spaces
are defined and studied in  Section
\ref{sectrect}.
To our knowledge, this is the first systematic study of this class of  sets, 
whose definition is inspired by some recent works of Lud\v ek Zaj\' \i \v cek
in \cite{Z}.
We extend the definition to Frechet spaces in Section \ref{small}
and prove in this more general setting that rectifiable sets 
of positive codimension are Baire-meager, Haar-null and Aronszajn-null.
The proof follows \cite{Z}.
In section \ref{comp}, we study the action of differentiable mappings
on rectifiable sets in separable Banach spaces.
We believe that this study is of independent interest, and hope
that it will have other applications.
The proof of Theorem \ref{main}  (and of a more general statement)
is then exposed in section
\ref{proof}. Actually, it consists mainly of stating appropriately
some of the ideas developed by Ma\~n\'e in
\cite{Mane:96} in order to reduce Theorem \ref{main}
to  an old result of Zaj\' \i \v cek
on monotone set-valued maps, see \cite{Z:78}.
This approach gives a more precise answer with an easier proof
than the finite dimensional approximation used in \cite{Mane:96}
and  in \cite{BC}.

\section{Rectifiable sets of finite codimension in
 Banach spaces}\label{sectrect}

In a finite dimensional Banach space $\Rm^n$,
one can say that a subset $A$ is countably rectifiable of codimension
$d$ if there exist countably many Lipschitz maps 
$F_i:\Rm^{n-d}\lto \Rm^n$
such that $A$ is contained in the union of the ranges of the maps 
$F_i$.
Many authors also add a set of zero $(n-d)$-dimensional
Hausdorff measure, but we do not.

In an infinite dimensional Banach space $B$, a first attempt
might be to define a rectifiable set of codimension $d$   as a set contained
in the countable union of ranges of Lipschitz
maps $F_i:B_i\lto B$, 
where $B_i$ are closed subspaces of $B$ of  codimension
$d$.
A closer look shows that this definition does not prevent 
$B$ itself from being rectifiable of positive codimension.
For instance, if $B$ is  a separable Hilbert space,
then $B\times \Rm^n$ is also a separable Hilbert space.
Therefore, it is isomorphic to $B$, and
 there exists a Lipschitz (linear) map  $B\lto B\times \Rm^n$ 
which is onto.
We thus  need  a finer definition, and the recent work of
Lud\v ek Zaj\' \i \v cek
in \cite{Z} opens the way.

A continuous linear map $L:B\lto B_1$ is called Fredholm
if its kernel is finite dimensional and if its range is closed
and has finite codimension.
We say that $L$ is a Fredholm linear map of type $(k,l)$
if $k$ is the dimension of the kernel of $L$
and $l$ is the codimension of its range. The index of $L$
is the integer $k-l$. Recall that the set of Fredholm linear 
maps is open in the space of continuous linear maps (for the norm
topology), and that the index is locally constant, although
the integers $k$ and $l$  are not.
They are lower semi-continuous.
To better understand the meaning of the index,
observe that, when $B$ and $B_1$ have finite dimension $n$ and $n_1$,
then the index of all linear maps is $i=n-n_1$.

\begin{defn}\label{rect}
Let $B$ be a Banach space.
We say that the  subset $A\subset B$ is a Lipschitz graph
of codimension $d$ if there exists:
\begin{itemize}
\item a splitting 
$B=D\oplus T$ of $B$, where $T$ is a linear subspace of dimension $d$
and $D$ is a closed linear subspace,
\item  a subset $D_1$ of $D$,
\item a Lipschitz map $g:D_1\lto T$,
\end{itemize}
such that
$$
A=\{(x_1\oplus g(x_1)): x_1\in D_1\}.
$$
We then say that $A$ is a Lipschitz graph transverse to $T$.

We say that the set $A\subset B$ is a rectifiable set
of codimension $d$ if there exist an integers $k$,
a Banach space $B_1$,
and a Fredholm
linear map $P:B_1\lto B$  of type $(k,0)$ such that
$$
A\subset P(A_1),
$$ 
where $A_1\subset B_1$ is a  Lipschitz graph of codimension $d+k$.

Finally, we say that the subset $A\subset B$ is countably
rectifiable of codimension $d$ if it is contained in the union
of countably many rectifiable sets of codimension $d$.
\end{defn}

This definition is directly inspired  by a recent work of Zaj\'i\v cek \cite{Z},
who proves that the  rectifiable sets of positive codimension according
to this definition are  small sets (see Section \ref{small} for more details).
 In particular, the space $B$  is not countably rectifiable of positive codimension 
 in itself.
This legitimates the systematic study of these sets initiated in the 
present paper.
Denoting by $\mR_d(B)$ the collection of all 
countably rectifiable subsets of codimension $d$ in $B$,
we have $\mR_{d+1}(B)\subset \mR_{d}(B)$.
This requires a proof:

\begin{lem}
If $A$ is countably rectifiable of codimension $d+1$,
then it is countably rectifiable of codimension $d$.
\end{lem}

\proof
We first prove that a Lipschitz graph $A$ of codimension 
$d+1$ is rectifiable of codimension $d$.
Let $T$ be a transversal to $A$, and let $D$ be a complement 
of $T$ in $B$.
Let $g:D\supset D_1\lto T$ be a Lipschitz map 
 such that $A=\{x\oplus g(x):x\in D_1\}$.
Let $S$ be a one dimensional subspace of $T$.
The Lipschitz graph
$A_2\subset D\times S\times T$ defined by
$$
A_2=\{(x,\lambda, g(x)-\lambda): x\in D_1, \lambda\in S\}
$$
has codimension $d+1$.
On the other hand, we have 
 $A\subset P(A_2)$, where $P:D\times S\times T\lto B$
is defined by
$$
P:(x,\lambda, t)\lmto x+\lambda+t
$$
which is of type $(1,0)$. As a consequence, $A$ is rectifiable of 
codimension $d$.

Assume now that $A$ is rectifiable of codimension $d+1$, and write 
it $A=L(\tilde A)$, where $\tilde A $ is a Lipschitz graph of
codimension $d+i+1$ and $L$ is linear Fredholm of type $(i,0)$.
Since $\tilde A$ is rectifiable of codimension $d+i$,
it can be written $\tilde A=P(\tilde A_1)$,
where $P$ is linear Fredholm of type $(l,0)$ and $\tilde A_1$
is a Lipschitz graph of codimension $d+l+i$.
Now we have $A=L\circ P(\tilde A_1)$, and $L\circ P$
is linear Fredholm of type $(i+l,0)$.
As a consequence, $A$ is rectifiable of codimension
$d$. 
\qed

Let us describe the action of Fredholm linear maps on rectifiable sets.

\begin{lem}\label{forward}
Let $L:B\lto B_1$ be a linear Fredholm map of index $i$,
and let $A\subset B$ be a rectifiable subset of codimension $d$.
Then $L(A)$ is rectifiable of codimension $d-i$ in $B_1$. 
\end{lem}

\proof
If $A=P(A')$, where $A'\subset B'$
is a Lipschitz graph of codimension $d+k$ and $P:B'\lto B$ is Fredholm
of type $(k,0)$, then $L(A)=L\circ P(A')$, and $L\circ P$ has index 
$k+i$. So it is enough to prove the statement when $A$ is a Lipschitz 
graph.

We now assume that $A$ is a Lipschitz graph of codimension $d$.
Let $K$ be the kernel of $L$, let $\breve K$ be a complement of $K$
in $B$, let $R$ be the range of $L$ and $\breve R$
be a complement of $R$ in $B_1$.
The set $A\times 0$ is a Lipschitz graph of codimension
$d+\dim \breve R$ in $B\times \breve R$.
On the other hand, the set $L(A)$ can also be written
$\tilde L(A\times 0)$, where 
$\tilde L:B\times \breve R\lto B_1$ is defined by
$\tilde L(b,r)=L(b)+r$.
The linear map $\tilde L$ is Fredholm of type 
$(\dim K,0)$, hence $L(A)=\tilde L(A\times 0)$
is rectifiable of codimension $d+\dim R-\dim K=d-i$.
\qed

\begin{lem}\label{preimage}
Let $L:B_1\lto B$ be a linear  map 
between two Banach spaces $B_1$ and $B$.
Assume that $\ker L$ has a closed  complement in $B_1$
and that the range of $L$ is closed and has finite codimension 
$l$.
If  $A$ is a countably rectifiable set of 
codimension $d$ in $B$, then $L^{-1}(A)$ is countably rectifiable
of codimension $d-l$ in $B_1$.
\end{lem}

\proof
Let $R\subset B$ be the range of $L$. The set $A\cap R$
is countably rectifiable of codimension $d$ in $B$.
Since $A\cap R=\pi(A\cap R)$, where
$\pi:B\lto R$ is a linear projection onto $R$,
and since $\pi$ is a Fredholm map of type $(l,0)$,
we conclude that $A\cap R$ is countably rectifiable
of codimension $d-l$ in $R$.
Let us now consider a splitting $B_1=R_1\oplus K_1$, where 
$K_1$ is the kernel of $L$.
Let $\tilde L:R_1\lto R$ be the restriction of $L$ to $R_1$.
Note that $\tilde L$ is a linear isomorphism, and therefore
$\tilde L^{-1}(A\cap R)$ is countably rectifiable of codimension
$d-l$  in $R_1$.
The conclusion now follows from
the observation that $L^{-1}(A)=\tilde L^{-1}(A\cap R)\times K_1.$
\qed

It is obvious from the definition that $\mR_d(B)$
is a translation-invariant $\sigma$-ideal of subsets of $B$.
More precisely, we have:
\begin{eqnarray*}
A\in \mR_d(B), A'\subset A 
&\Longrightarrow & \quad  A'\in \mR_d(B),
\\
A_n \in \mR_d(B)\;\forall n\in \Nm  &\Longrightarrow & \quad
\cup_{n\in \Nm} A_n \in \mR_d(B),
\\
A\in \mR_d(B),b\in B   &\Longrightarrow &\quad
b+A\in \mR_d(B).
\end{eqnarray*}
When $B=\Rm^n$ a countably rectifiable set of codimension
$d$ is what it should be:
a set which is contained in the union of the ranges
of countably many Lipschitz maps $f_i:\Rm^{n-d}\lto \Rm^n$.
Indeed such an range can be written as the projection 
on the second factor of the graph of $f_i$ in 
$\Rm^{n-d}\times \Rm^n$.
Some relations between finite-dimensional rectifiable sets
and infinite dimensional rectifiable sets of finite codimension
are given in \ref{presec}.
They are used in Section \ref{small} 
to prove, following Zaj\' \i \v cek
(see \cite{Z}), that rectifiable sets of finite codimension are small.
A consequence of these results is that a countably rectifiable 
set of positive codimension has empty interior.
As a consequence, we obtain:

\begin{lem}\label{subspace}
In a Banach space,
a closed linear subspace of codimension $d$ is rectifiable of codimension
$d$, but not countably rectifiable of codimension $d+1$.
\end{lem}

\proof
Let $B_1$ be a closed subspace of codimension $d$
in $B$. We can see $B_1$ as the range of a linear Fredholm
map $P:B\lto B_1$ of type $(d,0)$.
Since $B_1=P(B_1)$, if $B_1$ was countably rectifiable
of codimension $d+1$ in $B$, it would be countably
rectifiable of codimension $1$ in itself, which is 
in contradiction with the fact that countably rectifiable
sets of positive codimension have empty interior.
\qed

\subsection{Lipschitz graphs of finite codimension}
We  collect here some classical useful  facts concerning
Lipschitz graphs. Given a closed linear subspace $T$ of a Banach space $B$,
we consider the quotient space $B/T$ and the canonical projection 
$\pi:B\lto B/T$.
We endow  $B/T$ with the quotient norm $\|.\|$ defined by
$\|y\|:=\inf_{x\in \pi^{-1}(y)} \|x\|$.
It is well-known that $B/T$ is then itself a Banach space.
If $B_1$ is a closed complement of $T$, then $\pi_{|B_1}$ is a Banach space 
isomorphism onto $B/T$.

\begin{prop}
The following statements are equivalent for a subset $A\subset B$
and a finite dimensional (or more generally closed with a closed complement)
subspace $T$ of $B$:
\begin{enumerate}
\item
There exists a closed complement $D$ of $T$ in $B$,
 a subset $D_1\subset D$
and a Lipschitz map $g:D_1\lto T$ such that 
$A=\{x\oplus g(x):x\in D_1\}.
$
\item 
For each closed complement $D$ of $T$ in $B$, 
there exists a subset $D_1\subset D$
and a Lipschitz map $g:D_1\lto T$ such that 
$A=\{x\oplus g(x):x\in D_1\}.$
\item
The restriction to $A$ of the natural projection 
$\pi:B\lto B/T$ is a bi-Lipschitz homeomorphism onto its image.
\end{enumerate}
In this case, we say that $A$ is a Lipschitz graph transverse to $T$.
\end{prop}

\proof
If 1. holds, then $\pi(A)=\pi(D_1)$, and the 
restriction $\pi_{|A}$ can be inverted by the Lipschitz map
$$x\lmto \pi_{|D}^{-1}(x)\oplus g(\pi_{|D}^{-1}(x))
$$
so we have 3.

Let us now assume 3. If $D$ is a complement of $T$, then
we have 2. with 
$D_1=\pi_{|D}^{-1}(\pi(A))$
and 
$$
g=P\circ \pi_{|A}^{-1}\circ \pi_{|D},
$$
where $P:B\lto T$ is the projection associated to the splitting
$B=D\oplus T$. This map is Lipschitz since we have assumed that 
$\pi_{|A}^{-1}$ is Lipschitz.
\qed

\begin{prop}\label{delta}
Let $A$ be a Lipschitz graph transverse to $T$ in $B$.
Then there exists $\delta>0$ such that, if 
$F:B\lto B$ is Lipschitz with $\text{Lip}(F)<\delta$,
then $(Id+F)(A)$ is a Lipschitz graph  transverse to $T$.
\end{prop}

\proof
We consider a complement $D$ of $T$ in $B$, a subset $D_1$
of $D$ and a Lipschitz map $g:D_1\lto T$
such that $A=\{x\oplus g(x):x\in D_1\}$.
Let us set $G(x)=x\oplus g(x)$.
Let $\pi:B\lto B/T$ be the canonical projection.
We want to prove that $\pi\circ (Id+F)$
restricted to $A$ is a bi-Lipschitz homeomorphism.
It is equivalent to prove that 
$$\pi\circ(Id+F)\circ G=Id+ \pi\circ F\circ G$$
is a bi-Lipschitz homeomorphism of $D_1$  onto its image.
This holds if $\text{Lip}(\pi\circ F\circ G)<1$
by the classical inversion theorem for Lipschitz maps.
But $\text{Lip}(\pi\circ F\circ G)\leq \text{Lip}(F)\text{Lip}(G)$.
\qed

We have two useful corollaries:
\begin{lem}\label{diffeo}
Let $A\subset B$ be a Lipschitz graph
of codimension $d$, and let 
$F:U\lto B_1$ be a $C^1$ diffeomorphism onto its range,
where $U$ is an open subset of $B$.
For each $a\in A\cap U$, there exists an open neighborhood 
$V\subset U$ of $a$ such that $F(A\cap V)$
is a Lipschitz graph of codimension $d$.
\end{lem}

\begin{lem}\label{open}
Given a Lipschitz graph $A\subset B$ of codimension $n$, 
the set of tranversals to $A$ is open for the natural topology
on $n$-dimensional  linear subspaces of $B$.
\end{lem}

\textsc{Proof of Lemma \ref{open}:}
Let $\mG^d(B)$ be the set of $d$-dimensional linear subspaces of $B$.
Let $\mL(B)$ be the set of bounded linear selfmaps of $B$.
Let $A$ be a Lipschitz graph transverse to $T\in \mG^d(B)$.
By Proposition \ref{delta}, there exists $\delta>0$ such that 
$(Id+L)(A)$ is a Lipschitz graph transverse to $T$ when $L\in \mL(B)$
satisfies $\|L\|\leq \delta$.
We can assume that $\delta<1$, which implies that the linear map 
$(Id+L)$ is a Banach space isomorphism.
As a consequence, the set $A=(Id+L)^{-1}\big[(Id+L)(A)\big]$ is a Lipschitz
graph transverse to $(Id+L)^{-1}(T)$.
In other words, 
for each $L\in \mL(B)$ such that $\|L\|\leq \delta$, the space
$(Id+L)^{-1}(T)$ is a transversal to $A$.
These spaces, when $L$ vary, form a neighborhood of $T$ in $\mG^d(B)$
because the linear maps $\{(Id+L)^{-1},\|L\|\leq \delta\}$ form a neighborhood 
of the identity in $\mL(B)$.
\qed

\subsection{Pre-transversals of Rectifiable sets}\label{presec}
Let us now return to rectifiable sets which are not necessarily 
Lipschitz graphs.

\begin{defn}\label{pre}
Let $A$ be a rectifiable subset of codimension $d$ in the
Banach space $B$.
We say that the subspace $Q\subset B$ is a pre-transversal 
of $A$ if there exists:
\begin{itemize}
\item A Banach space $B_1$ and a Fredholm map $P:B_1\lto B$
of type $(k,0)$.
\item A Lipschitz graph $A_1$ such that $P(A_1)=A$.
\item A transversal $T$ of $A_1$ in $B_1$ such that 
$P_{|T}$ is an isomorphism onto $Q=P(T)$.
\end{itemize}
The dimension of  $T$ and $Q$ is necessarily  $d+k$.
\end{defn}

\begin{lem}\label{Q}
Let $A$ be a rectifiable subset of codimension $d\geq 1$ in the
Banach space $B$.
Then there exists an integer $n\geq d$
such that the  set of pre-transversals of $A$ contains 
 a non-empty open set
of $\mG^n(B)$ (the set of all $n$-dimensional linear subspaces of $B$).
\end{lem}

\proof
Let us write $A=P(A_1)$, where $P:B_1\lto B$
is a Fredholm map of type $(k,0)$  and $A_1\subset B_1$
is a Lipschitz graph of codimension $d+k$.
Let $K$ be the kernel of $P$, and 
let $T_0$ be a transversal of $A_1$ such that $T_0\cap K=0$,
such a transversal exists by Lemma \ref{open}.
Then, by definition, $P(T_0)=Q_0$ is a pre-transversal of 
$A$.
Let $B_2$ be a complement of $K$ in $B_1$ containing $T_0$,
and let  $U\subset \mG^{d+k}(B_2)$ be a neighborhood
of $T_0$ in the space of $(d+k)$-dimensional subspaces of $B_2$.
If $U$ is small enough, then each $T\in U$ (seen as a subspace of $B_1$) is a 
transversal of $A_1$
such that $T\cap K=0$, so that $Q=P(T)$ is a pre-transversal of $A$.
Since $P_{|B_2}$ is an isomorphism, the spaces $P(T), T\in U$
form a neighborhood of $Q_0$ in $\mG^{d+k}(B)$. 
\qed

\begin{lem}
Let $A$ be a rectifiable subset of codimension $d$ in the
Banach space $B$.
If $Q$ is a pre-transversal to $A$, then $(A+x)\cap Q$
is rectifiable of codimension $d$
(or equivalently it is rectifiable of dimension $(\dim Q)-d$) 
in the finite dimensional space $Q$
for each $x\in B$.
\end{lem}
\proof
We have $A=P(A_1)$ and $Q=P(T)$, where $P:B_1\lto B$
is Fredholm of type $(k,0)$ and where $A_1$ is a Lipschitz graph
transverse to $T$.
Let $K\subset B_1$ be the kernel of $P$, we have $K\cap T=0$.
The relation
$$
(A+x)\cap Q=P\Big((A_1+y)\cap (T\oplus K)\Big)
$$
holds for each point $y\in P^{-1}(x)$.
The set $(A_1+y)\cap (T\oplus K)$ is a Lipschitz graph of 
dimension at most $k=\dim K$ in $T\oplus K$.
As a consequence, the set 
$P((A_1+y)\cap (T\oplus K))$ is rectifiable of dimension 
at most $k$ in $Q$. In other words, it is rectifiable 
of codimension $d$ (recall that $\dim Q=d+k$).
\qed

The closure of a Lipschitz graph of codimension $d$ is
obviously a Lipschitz graph of codimension $d$.
This property does not hold for rectifiable sets,
but we have:

\begin{lem}\label{fsigma}
Each rectifiable set of codimension $d$ is contained
in a rectifiable set of codimension $d$ which is 
a countable union of closed sets.
\end{lem}

\proof
Using the notations of Definition \ref{rect}, we have 
$A=P(A_1)$, where $A_1$ is a Lipschitz graph of codimension 
$d$ in $B_1$.
The closure $\bar A_1$ of $A_1$ is a Lipschitz graph
of codimension $d$, hence $P(\bar A_1)$ is  rectifiable of
 codimension $d$.
Let us write $\bar A_1=\cap_{i\in \Nm} A^i_1$
where the sets $A^i_1, i\in \Nm$ are closed bounded 
subsets of $\bar A_1$.
Note that the sets $A^i_1$ are Lipschitz graphs.
We claim that $P(A^i_1)$ is closed for each $n$,
which implies the thesis.
Consider a sequence 
$a_n$ of  points of $P(A^i_1)$ converging to $a$ in $B$.
We want to prove that $a\in P(A^i_1)$.
Let $\tilde B_1$ be a closed complement of $\ker P$ in $B_1$, and
let $x_n$ be the sequence of preimages of $a_n$
in $\tilde B_1$. Since $P_{|\tilde B_1}$ is an isomorphism,
the sequence $x_n$ has a limit $x$ in $\tilde B_1$
with $P( x)=a$.
The point $x$ does not necessarily belong to $A_1^i$,
but  there exists a sequence $k_n$
in $\ker P$ such that $ x_n\oplus k_n\in A_1$ and such that 
$P(  x_n\oplus k_n)=a_n$.
The sequence $k_n$ is bounded because $A_1$ is bounded,
and therefore it has a
subsequence converging to a  limit $k$ 
in the finite dimensional space $\ker P$.
We have $ x\oplus k\in A^i_1$ because $A^i_1$ is closed, and 
thus $a=P( x\oplus k) $ belongs to $P(A^i_1)$.
 \qed

\section{Rectifiable sets in Frechet spaces}
\label{small}
In this section, we work in the more general setting of
Frechet spaces.
A Frechet space $F$ is a complete topological vector space
whose topology is generated by a countable family
of semi-norms $\|.\|_k, k\in \Nm$. 
Equivalently, the topology on $F$ is generated by
a translation invariant metric which turns $F$ into a complete
metric space. The main example here is $C^{\infty}(M)$.
Let us first define Lipschitz graphs and rectifiable sets 
in the context of Frechet spaces.

\begin{defn}\label{fgraphdef}
We say that the set $A\subset F$
is a Lipschitz graph of codimension $d$ if there exists 
a continuous linear map $P:F\lto B$
with dense range in the Banach space $B$
such that $P(A)$ is a Lipschitz graph of codimension $d$ in $B$.
\end{defn}

This definition is coherent with the already existing one
in the case where $F$ is a Banach space in view of the following:

\begin{prop}\label{fgraph}
Let $F$ and $B$ be Banach spaces, and let $P:F\lto B$
be a continuous linear map with dense range.
If $A\subset B$ is a Lipschitz graph  of codimension $d$ in $B$, then 
$P^{-1}(A)$ is a Lipschitz graph  of codimension $d$ in $F$.
\end{prop}

\proof
Let $T$ be a transversal of $A$ which belongs to the range of $P$
(such a transversal exists by Lemma \ref{open}),
and let $\tilde B$ be a complement of $T$ in $B$.
Let $T'\subset F$ be a subspace such that $P_{|T'}$ is an isomorphism
onto $T$, and let $\tilde F$ be the preimage of $\tilde B$.
Note  that $F=\tilde F\oplus T'$. 
This can be proved as follows:
let $\pi:B\lto T$ and $\tilde \pi:B\lto \tilde B$ be 
the projections corresponding to the splitting 
$B=T\oplus \tilde B$.
Each point $f\in F$ can be written 
$
f=t+(f-t)
$
with $t=P_{|T'}^{-1}\circ \pi \circ P (f)\in T'.$
It is enough to see that $P(f-t)\in \tilde B$.
This inclusion holds because
$P(f-t)=P(f)-\pi(P(f))=\tilde \pi(P(f))\in \tilde B$.
Let $g:\tilde B\supset D\lto T$ be the Lipschitz map
such that $A=\{x\oplus g(x), x\in D\}$.
Setting 
$D'=P^{-1}(D)$ and $g'=P^{-1}_{|T'}\circ g\circ P$,
we get that 
$$P^{-1}(A)=\{
f\oplus g'(f), f\in D'
\}
$$
and $g':D'\lto T'$ is Lipschitz.
\qed

\begin{prop}\label{frect}
The following statements are equivalent for a subset $A$ of the Frechet 
space $F$:
\begin{itemize}
\item
There exists a Banach space $B$ and a continuous linear map $P:F\lto B$
with dense range such that $P(A)$ is rectifiable of 
codimension $d$ in $B$. 
\item
There exists a Frechet space $F_1$,
a continuous linear map $\pi_1$ which is Fredholm of type $(k,0)$,
and a Lipschitz graph $A_1$ of codimension $d+k$ in $F_1$
such that $A=\pi_1(A_1)$.
\end{itemize}
We say that $A$ is rectifiable of codimension $d$ if these properties are 
satisfied.
\end{prop}

\proof
Assume that there exists a linear Fredholm map 
$\pi_1:F_1\lto F$ as in the statement.
Then, there exists a Banach space $B_1$ and a continuous
linear map $P_1:F_1\lto B_1$ with dense range such that 
$P_1(A_1)$ is a Lipschitz graph of codimension $d+k$ in 
$B_1$.
The space $K:=P_1(\ker \pi_1)\subset B_1$
has dimension at most $k$.
 Let us set $B:= B_1/K$, and let $\pi:B_1\lto B$
 be the standard projection.
 There is a unique map 
$P:F\lto B$ such that 
$P\circ \pi_1=\pi\circ P_1$.
In other words, the following diagram commutes.
$$
\begin{CD}
F_1 @>{P_1}>> B_1\\
@V{\pi_1}VV @V{\pi}VV\\
F @>{P}>> B
\end{CD}
$$
Let us check that $P$ has dense range and that 
$P(A)$ is rectifiable of codimension $d$.
The range of $P$ is the range of $\pi\circ P_1$,
which is dense because $\pi$ is onto and $P_1$ has dense range.
We have $P(A)=P(\pi_1(A_1))=\pi(P_1(A_1))$
which implies that $P(A)$ is rectifiable of codimension $d$
because $P_1(A_1)$ is a Lipschitz graph of codimension $d+k$
and  $\pi$ has type $(k',0)$ with $k'\leq k$.

Conversely, assume that $P:F\lto B$ exists as in the statement.
Then there exists a Lipschitz graph $A'_1$ of codimension $d+k$
in some Banach space $B_1$
 and a linear Fredholm map $\pi:B_1\lto B$
of type $(k,0)$ such that $P(A)=\pi(A'_1)$.
Let $K$ be the kernel of $\pi$, which has dimension $k$,
and let us set $F_1:=F\times K$.
Let $\pi_1:F_1\lto F$ be the projection on the first factor.
In order to complete the diagram with a map $P_1:F_1\lto B_1$,
we consider a right inverse $L$ of $\pi$,
(which is a Fredholm linear map of type $(0,k)$)
and set $P_1(f,k)=L(P(f))+k$.
Note then that $P\circ \pi_1=\pi\circ P_1$.
The range of $P_1$ is $\pi^{-1}(P(F))$, it is dense because
the range of $P$ is dense.
As a consequence, the set $A_1:=P_1^{-1}(A'_1)$
is a Lipschitz graph of codimension $d+k$, by definition.
On the other hand, the map $\pi_1$ is  Fredholm
of type $(k,0)$ hence 
the thesis follows
from the inclusion $A\subset \pi_1(A_1)$.
In order to prove this inclusion, let us consider a point $a\in A$.
There exist
$a'_1\in A'_1$ such that $\pi(a'_1)=P(a)$.
Then there exists $a_1\in A_1$ such that 
$\pi\circ P_1(a_1)=P(a)$,
which implies that $P(\pi_1(a_1))=P(a)$.
As a consequence, the difference $f:=a-\pi_1(a_1)$
belongs to $\ker P$.
Let us consider the point $b_1=a_1+(f,0)\in F_1$.
We have $P_1(b_1)=P_1(a_1)\in A'_1$, thus 
$b_1\in A_1$. On the other hand, 
$\pi_1(b_1)=\pi_1(a_1)+f=a$
hence $a\in \pi_1(A_1)$.
\qed

A subset $A\in F$ is said countably rectifiable
of codimension $d$ if it is a countable union of rectifiable
sets of codimension $d$.
The special case $F=C^{\infty}(M)$
may help to understand the definitions.

\begin{lem}
A subset $A\in C^{\infty}(M)$ is rectifiable of codimension 
$d$ if and only if there exists $p\in \Nm$
and a set $A'\subset C^p(M)$ which is rectifiable of codimension 
$d$ and such that $A=C^{\infty}(M)\cap A'$.
\end{lem}

\proof
Assume that $A$ is rectifiable of codimension $d$.
Then there exists a Banach space $B$ and a continuous 
linear map $P:C^{\infty}(M)\lto B$ with dense range such that 
$A'=P(A)$ is rectifiable of codimension $d$.
Then, the map $P$ is continuous for some $C^p$ norm and extends
to a continuous linear map $P_p:C^p(M)\lto B$ for some $p$.
Since the map $P_p$ has dense range, the set 
$A^p=P_p^{-1}(A')$ is rectifiable of codimension $d$,
and $A\subset A^p \cap C^{\infty}(M)$.

Conversely, if $A=C^{\infty}(M)\cap A'$
for some rectifiable set $A'\subset C^p(M)$,
then we have $A=P^{-1}(A')$, where $P:C^{\infty}(M)\lto C^p(M)$
is the standard inclusion. This inclusion is continuous with 
dense range hence $A$ is rectifiable.
\qed

We shall now explain, following 
Zaj\' \i \v cek
(see \cite{Z}) that rectifiable sets 
of positive codimension in Frechet spaces are small in various meanings of that
 term.
We first recall definitions.

A subset $A\subset F$ is called \textbf{Baire-meager}
if it is contained in a countable union of closed sets with empty interior.
Baire Theorem  states that a Baire-meager subset
of a Frechet space has empty interior.

A subset $A\subset F$ is called \textbf{Haar-null}
if there exists a compactly supported Borel  probability measure 
$\mu$ on $F$ such that $\mu(A+f)=0$ for all $f\in F$.
The equality $\mu(A+f)=0$ means that the set $A+f$
is contained in a Borel set $\tilde A_f$ such that 
$\mu(\tilde A_f)=0$.
A countable union of Haar-null sets is Haar-null, 
see \cite{Cr:72,BL} and \cite{prev} for the non-separable case.

A subset $A\subset F$ of a separable Frechet space $F$
is called \textbf{Aronszajn-null}
if, for each sequence $f_n$ generating a dense subset of $F$,
there exists a sequence $A_n$ of Borel subsets of $F$
such that $A\subset \cup_n A_n$ and such that,
for each $f\in F$ and for each $n$, the set
$$
\{x\in \Rm: f+xf_n\in A_n\}\subset \Rm
$$
has zero Lebesque measure. 
A countable union of Aronszajn-null sets is Aronszajn-null,
and each Aronszajn-null set is Haar null, see \cite{Ar:76,BL}.

\begin{thm}\label{Z}
Let $F$ be a  Frechet space, and let $A\subset F$
be countably rectifiable of positive codimension.  
Then  $A$ is  Baire-meager,  Haar null, 
and (if $F$ is separable)  Aronszajn-null.
\end{thm}

This Theorem is due to
 Lud\v ek Zaj\' \i \v cek
(see \cite{Z}) in the case of Banach spaces.
The extension to Frechet spaces that we now expose
is not very different.

\textsc{Proof of Theorem \ref{Z}:}
Let $A$ be a rectifiable set of positive codimension  in
the Frechet space  $F$, and let 
$P:F\lto B$ be a continuous linear map with dense image 
in the Banach space $B$ such that $P(A)$ is rectifiable 
of positive codimension.
We can assume without loss of generality that 
$A=P^{-1}(P(A))$.

Let us prove that $A$ is Baire meager.
By Lemma \ref{fsigma}, we can assume without loss 
of generality that $P(A)$ is a countable union
of closed sets in $B$, 
which implies that $A=P^{-1}(P(A))$ is a countable union of closed 
sets in $F$.
It is thus  enough to prove that $A$
has empty interior.
Let $\tilde Q$ be a pre-transversal of $P(A)$ in $B$  (see Definition \ref{pre}) 
contained in $P(F)$.
Such a space $\tilde Q$ exists by Lemma \ref{Q} because 
$P$ has dense range.
Let $Q\subset F$ be a linear subspace such that 
$P_{|Q}$ is an isomorphism onto $\tilde Q$.
Given $f\in F$, the set $(A+f)\cap Q$ has empty interior in $Q$,
and therefore there exists a sequence $q_n\in Q$
such that $q_n\lto 0$ (in $Q$ thus in $F$) and $(f+q_n)\not\in A$.
We conclude that the complement of $A$ is dense in $F$.

Let us prove that $A$ is Haar-null.
Let $\tilde Q$ be a pre-transversal of $A$ 
contained in $P(F)$ (see Definition \ref{pre} and Lemma \ref{Q})
and let $Q\subset F$ be a linear space such that 
$P_{|Q}$ is an isomorphism onto $\tilde Q$.
Let $\mu$ be the normalized Lebesgue measure on a bounded 
open subset of $Q$, that we also see as a compactly supported
Borel probability measure on $F$.
Since all rectifiable sets of positive codimension in the 
finite-dimensional space $Q$ have zero Lebesgue measure, we
conclude that  $\mu(A+f)=0$ for each $f\in F$. Therefore,  $A$ is 
Haar-null in $F$.

Finally, we assume that $F$ is separable and prove that 
$A$ is Aronszajn-null.
Let $f_n\in F$ be a sequence generating a dense subspace of $F$.
Since $P$ has a dense range, the sequence $P(f_n)$
generates a dense subset in $B$.
As a consequence, there exists an integer $N$
such that the space 
$$
\text{vect}\{P(f_n),n\leq N\}\subset B
$$ 
contains 
a pre-transversal $\tilde Q$ of $P(A)$.
Then, the space 
$$
F_N:=\text{vect}\{f_n,n\leq N\}\subset F
$$ 
contains a subspace $Q$ such that $P_{|Q}$ is an isomorphism
onto $\tilde Q$.
Since $(A+f)\cap Q$ has zero Lebesgue measure in $Q$
for each $f\in F$, we conclude from Fubini Theorem that
$
A\cap 
(F_n-f)
$
has zero Lebesque measure in $F_N$ 
 for each $f\in F$.
By standard arguments, (see for example
\cite{BL}, Proposition 6.29, p 144 or \cite{Ar:76}, 
Proposition 1, p 151)
 this implies that 
$A$ can be written as the union 
$$A=\cup_{n\leq N} A_n
$$
of Borel sets $A_n$ which are such that the set 
$$
\{ x \in \Rm: f+x f_n\in A_n\}
$$
has zero measure in $\Rm$ for each $f\in F$ and each $n\leq N$.
\qed

\section{Countably rectifiable sets and differential calculus}\label{comp}

In separable Banach spaces, the concepts of countably
rectifiable sets can be localized and behaves well 
with differential calculus. This 
section is not used in the proof of our results in Lagrangian dynamics, except a 
small part of it for Corollary \ref{epsilon}.

We say that a set $A\subset B$ is locally countably rectifiable
of codimension $d$ if, for each $a\in A$,
there exists a neighborhood $U$ of $a$ in $B$
such that $U\cap A$ is countably rectifiable of codimension $d$.
Every countably rectifiable set is obviously locally 
 countably rectifiable, and conversely we have:

\begin{lem}\label{local}
Let $B$ be a separable Banach space.
If the subset $A\subset B$ is locally countably rectifiable of 
codimension $d$, then it is countably rectifiable of codimension 
$d$.
\end{lem}
\proof
The space $A$ is a separable metric space for the metric induced
from the norm of $B$, and therefore it has the Lindel\"of
property: each open cover of $A$ admits a countable subcover.
Now the hypothesis of local countably rectifiability
implies that $A$ can be covered by a union of open subsets
of $A$ each of which is countably
rectifiable of codimension $d$ in $B$.
 Therefore, by the Lindel\"of property, $A$ is the  union
of countably many sets each of which is countably
rectifiable of codimension $d$ in $B$.
\qed

Let $B$ and $B_1$ be two Banach spaces, $U\subset B$
be an open subset of $B$, and $F:U\lto B_1$ be a $C^1$ map.
We say that $F$ is Fredholm if the Frechet differential
$dF_x$ is Fredholm at each $x\in U$.
If $U$ is connected, then the index of $dF_x$ does not depend
on $x$, we say that this is the index of $F$.
The following result shows that rectifiable sets of codimension $d$
in separable Banach spaces
could have been equivalently defined as the image by a $C^1$
Fredholm map of index $i$ of a Lipschitz graph of codimension 
$d+i$.

\begin{prop}\label{Fredholm1}
Let $B$ and $B_1$ be  separable Banach spaces,
 let $U$ be an open subset
of $B$, and let $A\subset U$ be a countably rectifiable set 
of codimension $d$.
If $F:U\lto B_1$ is a $C^1$ Fredholm map of index $i$
then $F(A)$ 
is countably rectifiable  of codimension $d-i$
in $B_1$.
\end{prop}

\proof
It is enough to prove that, for each $i$ and $d$,
 the image of a Lipschitz 
graph of codimension $d$ by a Fredholm map of index $i$
is locally countably rectifiable of codimension $d-i$.
In fact, if $A=\cup_n P_n(A_n)$ is a countably rectifiable set
of codimension $d$, where $A_n$ are Lipschitz graphs of codimension
$d+i_n$ and $P_n$ are linear Fredholm maps of type $(i_n,0)$,
then $F(A)=\cup _n F\circ P_n(A_n)$.
The map $F\circ P_n$ is Fredholm of index $i+i_n$, and therefore,
if our claim is proved, then 
$F\circ P_n(A_n)$ is countably rectifiable of codimension 
$(d+i_n)-(i+i_n)=d-i$.

So we now assume that $A$ is a Lipschitz graph.
In order to prove that the image $F(A)$ is countably rectifiable of codimension 
$d-i$, it is enough to prove that each point $a\in A$ has a neighborhood
$U$ in $A$ such that $F(U)$ is countably rectifiable of codimension $d-i$.
Let $a\in A$ be given. Assume that the linear Fredholm map
$dF_a$ is of type $(k,l)$. There is a local $C^1$
diffeomorphism $\phi$ of $B$ around $a$
and a splitting $B= B_1\oplus C$ of $B$,
with $\dim C=l$, 
such that 
$F=\tilde F\circ \phi$ in a neighborhood of $a$,
where $\tilde F$ is of the form
$$x\lmto \pi(x)\oplus f(x),
$$
where $\pi:B\lto B_1$ is a linear Fredholm map of type
$(k,0)$ and $f:B\lto C$ is a $C^1$ map satisfying 
$df_a=0$, see \cite{BZS}, Theorem 1.1.
By Lemma \ref{diffeo}, 
$\phi(A)$ is locally a Lipschitz graph of codimension $d$.
In other words, there exists a neighborhood $U$ of $a$ in $A$ such that 
$\phi(U)$ is a Lipschitz graph of codimension $d$.
Hence there exists a 
splitting $B=\tilde B\oplus T$, a subset 
$\tilde D\subset \tilde B$,
and a Lipschitz map $g:\tilde D\lto T$ such that 
$\phi(U)= \{x\oplus g(x): x\in \tilde D\}$.
Let us define
$ A_2\subset  B\times  C$
as
$$
A_2:= \big\{\big(x,f(x)\big): x\in \phi(U)\big\}
=\big\{\big(\tilde x\oplus g(\tilde x),f(\tilde x\oplus g(\tilde x))\big):
\tilde x\in \tilde D\big\}.
$$
Since the map 
$$
\tilde D\ni\tilde x\lmto (g(\tilde x),f(\tilde x\oplus g(\tilde x))
\in T\times C$$
is Lipschitz,
the set $A_2$
is a Lipschitz graph transverse to $T\times C$ in $B\times C$,
it is thus of codimension $d+l$.
Now the set $\tilde F(A)$ is just the image of $A_2$
by the linear map 
$$B\times C \ni (x,z)\lmto \pi(x)\oplus z \in B.
$$
This linear map is Fredholm of type $(k,0)$, so
 by definition $\tilde F(A)$
 is rectifiable of codimension
$d+l-k=d-i$.
\qed

The following direct corollary is especially important:
\begin{cor}\label{inv}
Let $B$ be a separable Banach space, $U\subset B$ and open set,
and $\phi:U\lto B$ a $C^1$ diffeomorphism onto its image $V$.
Then the set $A\subset U$ is countably rectifiable of codimension $d$
if and only if its image $\phi(A)$ is countably rectifiable of codimension $d$.
\end{cor}

A consequence of all these observations is that 
the notion of countably rectifiable subsets of 
codimension $d$ is well-defined in separable manifolds
modeled on separable Banach spaces:

\begin{defn}
Let $W$ be a separable $C^1$ manifold modeled on the separable
Banach space $B$.
The set $A\subset B$ is said countably rectifiable of codimension $d$ if, for each
 point
$a\in A$, there exists a neighborhood $U\subset W$ of a in $W$ and a chart 
$\phi:U\lto V\subset B$ such that $\phi(U\cap A)$ is countably rectifiable
of codimension $d$ in $B$.
\end{defn}

In view of Corollary \ref{inv}, if $A$ is a countably rectifiable set 
of codimension $d$ in $W$ and $\phi:U\subset W\lto V\subset B$ is any chart of $W$,
then $\phi(A\cap U)$ is countably rectifiable of codimension $d$ in $B$.
Obviously, when $W=B$, this definition of countably rectifiable sets of
 codimension $d$
coincides with the former one.
Proposition \ref{Fredholm1} has a straightforward generalization:

\begin{thm}\label{Fredholm}
Let $W$ and $W_1$ be  separable manifolds modeled on separable Banach spaces,
  and let $A\subset W$ be a countably rectifiable set 
of codimension $d$.
If $F:W\lto W_1$ is a $C^1$ Fredholm map of index $i$
then $F(A)$ 
is countably rectifiable  of codimension $d-i$
in $W_1$.
\end{thm}

If $W$ and $W_1$ are two separable manifolds modeled on separable Banach
 spaces $B$ and $B_1$, 
the   $C^1$ map $F:W\lto W_1$ is called a submersion
at $x$
if the differential $dF_x:T_xW\lto T_{F(x)}W_1$
 is onto, and if its kernel splits (if these conditions
 are satisfied, we say that $dF_x$ is a linear submersion).
It means that there exists a closed linear subspace 
$\tilde B$ in $T_xM$ such that $T_xM=\ker(dF_x)\oplus \tilde B$.
It is known that the map $F$ is a submersion at $x$
if and only if there exist local charts at $x$ and $F(x)$
such that the expression of $F$ in these charts is a linear 
submersion. The following statement then follows from Lemma 
\ref{preimage}:

\begin{prop}\label{submersion}
Let $W$ and $W_1$ be two separable manifolds
modeled on separable Banach spaces.
Let $A\subset W_1$ be a countably rectifiable 
subset of codimension $d$, and let $F:W\lto W_1$
be a $C^1$ map 
which is a submersion at each point of $F^{-1}(A)$.
Then $F^{-1}(A)$ is countably rectifiable of 
codimension $d$.
\end{prop}

\section{Application to Lagrangian systems}\label{proof}
We now return to the study of Minimizing measures of Lagrangian systems.
Let us begin with  a general abstract  result:
\begin{thm}\label{abstract}
Let $F$ be a  Frechet space of $C^2$ functions
on $TM$ (but not necessarily with the $C^2$ topology) and
$U$ be an open subset of $F$.
Assume that
\begin{itemize}
\item The topology on $F$ is stronger than the compact-open topology.
 In other words, for each compact set $K\subset TM$
the natural map $F\lto C(K)$ is continuous.
\item The space $F$ contains a dense subset of $C(M)$ (
where the functions of $C(M)$ are seen as functions on $TM$).  
\item For each $f\in U$, the function $L-f$
is a Tonelli Lagrangian.
\end{itemize}
Then, for each $k\in \Nm$, the set
$$
\{ f\in U : \dim(\mM(L-f))\geq k
\}
$$
is a countable union of Lipschitz graphs of codimension
$k$ in $F$ (and therefore it is countably rectifiable of codimension 
$k$).
\end{thm}

Before we turn to the proof,
let us  see how to derive the statements of the introduction
from this result.

By taking $U=F=C^p(M), p\in \{2,3,\ldots, \infty\}$,
we  obtain that the set $\sigma ^p(L)$
of functions $f\in C^p(M)$ such that $L-f$ has more than 
one Mather mesure
is countably rectifiable of codimension 1,
 which is stronger than  the Theorem of Ma\~n\'e 
 (Theorem \ref{Mane}).

\textsc{Proof of Theorem \ref{main}: }
As earlier, let us identify 
the space $H^1(M,\Rm)$ with a $b_1$-dimensional
space of smooth forms, and therefore with a 
$b_1$-dimensional space of funtions on $TM$.
Let us take 
$U=F=H^1(M,\Rm)\times  C^p(M)$, and apply Theorem
\ref{abstract}.
We obtain that the set of pairs 
$(c,f)\in H^1(M,\Rm)\times C^p(M)$ such that 
$\dim (\mM(L-c-f))\geq k$ is countably rectifiable of codimension
$k$. 
Since $\Sigma^p_k $ is the projection
of this set on the second factor $C^p(M)$, we conclude by Proposition 
\ref{frect} that $\Sigma^p_k$ is countably rectifiable of
codimension $k-b_1$.
\qed

\textsc{Proof of Corollary \ref{epsilon}: }
Let $P_p:C^p(M)\lto C^2(M)$ be the standard inclusion,
which is a continuous linear map with dense range.
We observe that 
$$
A_k^p(M)=A_k^2(M)\cap C^p(M)=P_p^{-1}(M)
$$ 
so, by Proposition \ref{frect},
it is enough to prove the result for $p=2$.
The map 
$$(0,\infty)\times C^2(M,\Rm)\ni (\epsilon,V)\lmto
\epsilon V\in C^2(M,\Rm)
$$
is a smooth submersion.
By theorem \ref{main}, the set 
$$\Sigma^2_k:=\{V\in C^2(M,\Rm) : \max_{c\in H^1(M,\Rm)}
\dim \mM (L-V+c)\geq k
\}
$$  
 is countably rectifiable of
codimension $k-b_1$, hence 
Proposition \ref{submersion} implies that the set 
of pairs $(\epsilon,V)\in (0,\infty)\times C^2(M,\Rm)$
such that $\epsilon V\in \Sigma^2_k$ is countably rectifiable of 
codimension $k-b_1$ in $(0,\infty)\times C^2(M,\Rm)$.
The conclusion now follows by Lemma \ref{forward}
since  $A^2_k$ is the projection of this set on the second 
factor.
\qed

The proof of  Theorem \ref{abstract}
will occupy the rest of the section.
It is useful first to recall the notion of a 
\textbf{closed
probability measure} on $TM$.
The Borel probability measure $\mu$ is said closed
if it is compactly supported and if,
 for each  function $f\in C^{\infty}(M)$, we have
$$\int_{TM}df_x\cdot v \;d\mu(x,v)=0.
$$
If $\mu$ is compactly supported and invariant under 
the Euler-Lagrange flow
of a Tonelli Lagangian $L$, then $\mu$ is closed.
We recall the proof first given by Mather, see \cite{Mather:91}.

We can express the invariance of the measure $\mu$
by saying that, for each smooth function $g:TM\lto \Rm$,
we have
$$\int_{TM}gd\mu=\int_{TM}g\circ \psi^t d\mu.
$$
Differentiating at $t=0$, we get 
$$\int_{TM} dg_{(x,v)}\cdot \mX(x,v) d\mu(x,v)=0,
$$
where $\mX$ is the Euler-Lagrange vector-field (the generator of 
$\psi^t$).
This is true in particular if $g=f\circ \pi$ is a function which
depends only on the position $x$, and writes
$$\int_{TM}df_x\cdot \Pi(\mX(x,v)) d\mu(x,v)=0
$$
for each smooth function $f$ on $M$,
where $\Pi:T_{(x,v)}(TM)\lto T_x M$ is the 
differential of the standard projection.
The proof now follows from the observations that $\Pi(\mX(x,v))=v$.

Let us now fix a Riemannian metric on $M$,
and define, for each $n\in \Nm$, the compact subset 
$K_n$ of $TM$ as follows:
$$K_n=\{(q,v)\in TM : \|v\|_q\leq n\}.
$$ 
We denote by $\mC_n$ the set of closed probability measures
supported on $K_n$.
We also define the Banach space $B_n$
as the closure, in $C(K_n)$, of the restrictions 
to $K_n$ of the functions of $F$.
Since $C(K_n)$ is separable, so is $B_n$.
Moreover, it follows from our assumptions on $F$
that $B_n$ contains $C(M)$.
Each probability measure $\mu$ on $K_n$
gives rise to a linear form 
on $C(K_n)$ and, by restriction,
to a linear form $l_n(\mu)\in B_n^*$,
 which is just defined by
$$
l_n(\mu)\cdot f=\int fd\mu.
$$
Because  $B_n$ is not necessarily dense in $C(K_n)$,
the map $l_n$ is not necessarily one-to-one.

Let us now define the functional 
$A_n:B_n\lto \Rm$ by
$$
A_n(f)=\sup_{\mu\in \mC_n} \int (f-L)d\mu.
$$
This functional is convex, and it is bounded from below
(because any Dirac supported on a point of the zero section 
of $TM$ belongs to $\mC_n$).
By standard results of convex analysis, 
the set 
$$
\{f\in B_n : \dim(\partial A_n(f))\geq k
\}
$$
is a countable union of Lipschitz graphs of codimension 
$k$ in the separable Banach space  $B_n$ for each $k\in \Nm$
 (see \cite{Z:78} or \cite{BL}, Theorem 4.20, p. 93),
 where $\partial A_n (f)$ is the sub-differential 
 in the sense of convex analysis of the function  $A_n$
 at point $f$.
We will now be able to conclude if we can relate 
the set $\mM(L-f)$ of Mather measures and the set 
$\partial A_n(f)$.
This is the content of the following,
where we denote by $P_n :F\lto B_n$ the natural map :

\begin{lem}\label{n}
If $L-f$ is a Tonelli Lagrangian,
then there exists $n\in \Nm$ such that 
$$
\dim (\mM(L-f))\leq \dim \partial A_n(P_n(f)).
$$
\end{lem}

Assuming Lemma \ref{n}, let us finish the proof 
of Theorem \ref{abstract}.
The map $P_n:F\lto B_n$
is continuous by our assumptions on $F$,
and it has a dense range by the definition of $B_n$.
In view of the Lemma
we have 
\begin{equation}\tag{U}\label{U}
\{f\in U: \dim \mM(L-f)\geq k
\}\subset
\bigcup_n
P^{-1}_n\big(
\{f\in B_n: \dim \partial A_n(f)\geq k
\}\big).
\end{equation}
Each of the sets
$$\{f\in B_n: \dim \partial A_n(f)\geq k
\}
$$
is a  countable union of Lipschitz graphs of codimension $k$
in $B_n$,
hence the preimage 
$$
P^{-1}_n\big(
\{f\in B_n: \dim \partial A_n(f)\geq k
\}\big)
$$
is a countable union of Lipschitz graphs
of codimension $k$ in $F$, by Definition \ref{fgraphdef}.
Therefore Theorem \ref{abstract} follows from (\ref{U}).
\qed

The last step is to prove Lemma \ref{n}.
The main tool is 
 the following beautiful
variational principle which has been established
by Bangert \cite{Ba} and Fathi and Siconolfi \cite{FS:04}
following fundamental ideas of Ma\~n\'e,  \cite{Mane:96}
(see also \cite{indiana,ENS}):

\begin{thm}
The Mather measures of the Tonelli Lagrangian $L$
are those which minimize the action $\int Ld\mu$
in the class of all compactly supported closed measures.
\end{thm}

The key point here is that the invariance of the measure
is obtained as a consequence of its minimization property,
and not imposed as a constraint as in Mather's work.
By applying this general result to the Lagrangian
$L-f$, we obtain that there exists $n\in \Nm$
such that
the set $\mM(L-f)$ of Mather measures
is the set of measures $\mu\in \mC_n$
which minimize the action $\int (L-f)d\mu$ in $\mC_n$.
Now if $\mu$ is such a measure, then the 
associated  
$l_n(\mu)$ belongs to $\partial A_n(P_n(f))$, as can easily be seen
from the definition of $A_n$.
In other words, we have proved that 
$$
l_n(\mM(L-f))\subset \partial A_n(P_n(f))
$$
when $n$ is large enough.
In order to prove Lemma \ref{n}, it is enough 
to observe that $l_n$ is one to one on $\mM(L-f)$.
This property holds because $B_n$ contains $C(M)$
and because the elements of $\mM(L-f)$ are all supported on a 
Lipschitz graph.
\qed

\bibliographystyle{amsplain}
\providecommand{\bysame}{\leavevmode\hbox
to3em{\hrulefill}\thinspace}

\end{document}